\documentclass[11pt]{amsart}
\usepackage{amssymb}
\usepackage{amsfonts}
\usepackage{color}

\input{diagrams}
\diagramstyle[noPostScript]

\newarrow{Dash}{}{dash}{}{dash}>%....>
\newarrow{Into}C--->
\newarrow{Onto}----{>>}
\newarrow{Equal}===={=}
\newtheorem{theorem}{Theorem}[section]
\newtheorem{proposition}[theorem]{Proposition}
\newtheorem{lemma}[theorem]{Lemma}
\newtheorem{corollary}[theorem]{Corollary}
\newtheorem{definition}[theorem]{Definition}

\theoremstyle{remark}

\renewenvironment{proof}{{\noindent\bf Proof.}}{\hfill $\Box$\par\vskip3mm}

\newcommand{\Aa}{\mathcal{A}}
\newcommand{\Bb}{\mathcal{B}}
\newcommand{\Cc}{\mathcal{C}}

\newcommand{\Ll}{\mathcal{L}}
\newcommand{\Mm}{\mathcal{M}}

\newcommand{\Ss}{\mathcal{S}}

\def\NN{{\mathbb N}}

\begin{document}
\title{The finite {\it Rat}-Splitting for Coalgebras}

\begin{abstract}
Let $C$ be a coalgebra. We investigate the problem of when the rational part of every finitely generated $C^*$-module $M$ is a direct summand $M$. We show that such a coalgebra must have at most countable dimension, $C$ must be artinian as right $C^*$-module and injective as left $C^*$-module. Also in this case $C^*$ is a left Noetherian ring. Following the classic example of the divided power coalgebra where this property holds, we investigate a more general type of coalgebras, the chain coalgebras, which are coalgebras whose lattice of left (or equivalently, right, two-sided) coideals form a chain. We show that this is a left-right symmetric concept and that these coalgebras have the above stated splitting property. Moreover, we show that this type of coalgebras are the only infinite dimensional colocal coalgebras for which the rational part of every finitely generated left $C^*$-module $M$ splits off in $M$, so this property is also left-right symmetric and characterizes the chain coalgebras among the colocal coalgebras.
\end{abstract}

\author{Miodrag Cristian Iovanov}%$^*$}
\thanks{$^*$ This paper was partially supported by a CNCSIS BD-type grant and was written within the frame of the bilateral Flemish-Romanian project "New Techniques in Hopf Algebra Theory and Graded Rings"}
%\thanks{$^2$}
\date{}
\maketitle

\section*{Introduction}
%\noindent

Let $R$ be a ring and $T$ be a torsion preradical on the category of left $R$-modules ${}_{R}\Mm$. Then $R$ is said to have splitting property provided that $T(M)$, the torsion submodule of $M$, is a direct summand of $M$ for any $M\in{}_{R}\Mm$. More generally, if $\Cc$ is a Grothendieck category and $\Aa$ is a subcategory of $\Cc$, then $\Aa$ is called closed if it is closed under subobjects, quotient objects and direct sums. To every such subcategory we can associate a preradical $t$ (also called torsion functor) if for every $M\in \Cc$ we denote by $t(M)$ the sum of all subobjects of $M$ that belong to $\Aa$. We say that $\Cc$ has the splitting property with respect to $\Aa$ if $t(M)$ is a direct summand of $M$ for all $M\in\Cc$. In the case of the category of $R$-modules, the splitting property with respect to some closed subcategory is a classical problem which has been considered by many authors. In particular, when $R$ is a commutative ring, the question of when the (classical) torsion part of an $R$ module splits off is a well known problem. J. Rotman has shown in \cite{Rot} that for a commutative domain the torsion submodule splits off in every $R$-module if and only if $R$ is a field. I. Kaplansky proved in \cite{K1}, \cite{K2} that for a commutative integral domain $R$ the torsion part of every finitely generated $R$-module $M$ splits in $M$ if and only if $R$ is a Pr$\rm\ddot{u}$fer domain. While complete results have been obtained for commutative rings, the problem still remains wide open for the non-commutative case. %While complete results have been obtained in the commutative case, the characterization of the noncommutative rings $R$ for which the torsion splits in every $R$ module (or in every finitely generated module) is still an open problem.\\
In this paper we investigate the situation when the ring $R$ arises as the dual algebra of a $K$-coalgebra $C$, $R=C^*$. Then the category of the left $R$-modules naturally contains the category $\Mm^C$ of all right $C$-comodules as a full subcategory. In fact, $\Mm^C$ identifies with the subcategory $Rat({}_{C^*}\Mm)$ of all rational left $C^*$-modules, which is generally a closed subcategory of ${}_{C^*}\Mm$. Then two questions regarding the splitting property with respect to $Rat({}_{C^*}\Mm)$ naturally arise: first when is the rational part of every left $C^*$-module $M$ a direct summand of $M$ and when does the rational part of every finitely generated $C^*$-module $M$ split in $M$. The first problem, the splitting of ${}_{C^*}\Mm$ with respect with the closed subcategory $Rat({}_{C^*}\Mm)$ has been treated by C. N\u as\u asescu and B. Torrecillas in \cite{NT} where it is proved that if all $C^*$-modules split with respect to $Rat$ then the coalgebra $C$ must be finite dimensional. The techniques used involve some amount of category theory (localization in categories) and strongly relies on some general results of M.L.Teply from \cite{T1}, \cite{T2}, \cite{T3}. \\
We consider the more general problem of when $C$ has the splitting property only for finitely generated modules, that is, the problem of when is the rational part $Rat(M)$ of $M$ a direct summand in $M$ for all finitely generated left $C^*$-modules $M$. We call these coalgebras left finite $Rat$-splitting coalgebras (or we say that they have the left finite $Rat$-splitting property). If the coalgebra $C$ is finite dimensional, then every left $C^*$-module is rational so $\Mm^C$ is equivalent to ${}_{C^*}\Mm$ and $Rat(M)=M$ for all $C^*$-modules $M$ and in this case $Rat(M)$ trivially splits in any $C^*$-module. Therefore we will deal with infinite dimensional coalgebras, as generally the infinite dimensional coalgebras produce examples essentially different from the ones in algebra theory. We first prove some general properties for left finite $Rat$-splitting coalgebras, namely such a coalgebra $C$ is artinian as right $C^*$-module and injective as left $C^*$-module, it has at most countable dimension and has finite dimensional coradical. Also $C^*$ is a left Noetherian ring. We look at a very simple example of a coalgebra where this property holds, namely, the divided power coalgebra (see \cite{DNR}, Example 1.1.4), which has $K[[X]]$ as its dual algebra. This is in some sense the simplest possible example of infinite dimensional coalgebra that has the left (and right) finite $Rat$-splitting property. We introduce and study (left) chain coalgebras to be the coalgebras for which every two left subcomodules $M,N$ satisfy either $M\subseteq N$ or $N\subseteq M$. We see that this is a left-right symmetric concept and we give a simple characterization of these coalgebras as being exactly those having each factor of the coradical series a simple comodule. Moreover, this gives a complete characterization of these coalgebras in the case when the base field is algebraically closed: the divided power coagebra and its subcoalgebras are the only ones of this type. We show that chain coalgebras have the (left and right) finite $Rat$-splitting property. Then we investigate the colocal finite $Rat$-splitting coalgebras. In the main result of the paper we show that a colocal coalgebra satisfying the left finite $Rat$-splitting property must be a chain coalgebra, and therefore it also has the right finite $Rat$-splitting property. This provides a characterization of the divided power coalgebra over an algebraically closed field (or more generally of chain coalgebras) among local coalgebras, namely they are exactly those coalgebras $C$ for which the rational part of every finitely generated left (or right) $C^*$-module splits off. 

%In fact, the main result of the paper characterizes the infinite dimensional co-local coalgebras that satisfy the left finite splitting property; we prove that these are exactly the chain coalgebras, thus showing that thay also satisfy the right finite splitting property.
%\end{section}

\section{Splitting Problem}

Let $C$ be a coalgebra with counit $\varepsilon$ and comultimplication $\Delta$. We use the Sweedler convention $\Delta(c)=c_1\otimes c_2$ where we omit the summation symbol. For a vector space $V$ and a subspace $W$ of $V$ denote by $W^{\perp}=\{f\in V^*\mid f(x)=0,\,\forall\,x\in W\}$ and for a subspace $X\in V^*$ denote by $X^\perp=\{x\in V\mid f(x)=0,\,\forall\,f\in X\}$. If $M$ is a right (or left) $R$-module denote by $\Ll_R(M)$ (or ${}_{R}\Ll(M)$) the lattice of the submodules of $M$. Also, if $S$ is another ring and $Q$ is a fixed $R$-$S$-bimodule, for any left $R$-module $M$ we have applications 
$${}_{R}\Ll(M)\ni N\rightarrow N^\perp=\{f\in {\rm Hom}_{R}(M,Q)\mid f(x)=0,\,\forall\,x\in N\}\in \Ll_S({\rm Hom(M,Q)})$$
$$\Ll_S({\rm Hom(M,Q)})\ni X\rightarrow X^\perp=\{x\in M\mid f(x)=0,\,\forall\,f\in X\}\in {}_{R}\Ll(M)$$
forming a Galois pair (see \cite{AN}). 

\begin{lemma}\label{1.perp}
Let $C$ be a coalgebra. Then for any finitely generated right (or left, or two-sided) submodule $X$ of $C^*$, $(X^{\perp})^\perp=X$.
\end{lemma}
\begin{proof}  %f(c*d*)(c)=(c*d*)(c1)c2=c*(c1)d*(c2)c3; f(c*)f(d*)(c)=f(c*)(d*(c1)c2)=d*(c1)c*(c2)c3
%[THIS PROBABLY FOLLOWS ALSO FROM THE BOOK {\it Finiteness conditions in module Theory} ALBU, NASTASESCU]\\
Put $R={\rm End}(C^C,C^C)$. Then $R$ is a ring with multiplication "$\cdot$" equal to opposite composition of morphisms. Let $M={}_{C^*}C$ and $Q={}_{C^*}C_R$ where the right $R$-module structure on $C$ is $c\cdot f=f(c)$. It is not difficult to see that the isomorphism of rings $C^*\simeq {\rm End}(C^C,C^C)=R$, $c^*\rightarrow (c\mapsto c^*(c_1)c_2)$, transposes the problem to the Galois correspondence $X\rightarrow X^\perp$ between the left $R$ module $C$ and the right $R$ module ${\rm End}^C(C,C)={\rm Hom}({}_{C^*}C,{}_{C^*}C)$. That is, it is enough to prove the statement for finitely generated right ideals of $R$. Suppose $X\subseteq R$ is a right ideal generated by $f_1,\dots,f_n$ as right $R$ module and let $f\in R$ such that $f\vert_{X^\perp}=0$. Then we have $X=\sum\limits_{i=1}^{n}f_i\cdot R$ so $X^\perp=\bigcap\limits_{i=1}^{n}(f_i\cdot R)^\perp=\bigcap\limits_{i=1}^{n}{\rm Ker}f_i$. Then $f$ induces a morphism $\overline{f}:\frac{C}{\bigcap\limits_{i=1}^{n}{\rm Ker}f_i}\rightarrow C$ and as $C^C$ is injective, the canonic monomorphism  $0\rightarrow \frac{C}{\bigcap\limits_{i=1}^{n}{\rm Ker}f_i} \rightarrow \bigoplus\limits_{i=1}^{n}\frac{C}{{\rm Ker}f_i}$ gives rise to the exact sequence 
$$\bigoplus\limits_{i=1}^n{\rm Hom}^C(\frac{C}{{\rm Ker}f_i},C) \simeq  {\rm Hom}^C(\bigoplus\limits_{i=1}^n\frac{C}{{\rm Ker}f_i},C)\stackrel{\varphi}{\rightarrow}{\rm Hom}^C(\frac{C}{\bigcap\limits_{i=1}^{n}{\rm Ker}f_i}, C)\rightarrow 0$$
Let $(g_i)_{i={\overline{1,n}}}\in \bigoplus\limits_{i=1}^n{\rm Hom}^C(\frac{C}{{\rm Ker}f_i},C)$ be such that $\varphi(\sum\limits_{i=1,n}g_i)=\overline{f}$. As for any $i$ we have a monomorphism $\overline{f_i}:\frac{C}{{\rm Ker}f_i}\rightarrow C$ induced by $f_i$ and as the right $C$-comodule $C$ is injective, any the diagram
\begin{diagram}
0 & \rTo & \frac{C}{{\rm Ker}f_i} & \rTo^{\overline{f_i}} & C\\
  &      & \dTo^{g_i}             & \ldDotsto_{h_i}       &  \\
  &      & C                      &                       &  \\
\end{diagram}
can be completed commutatively by a morphism of right $C$-comodules $h_i$. Then we have $\varphi(\sum\limits_{i=1}^n\overline{f_i}\cdot h_i)=\varphi(\sum\limits_{i=1}^nh_i\circ\overline{f_i})=\overline{f}$ and composing this with the cannonical projection $p:C\rightarrow \frac{C}{\bigcap\limits_{i=1}^n{\rm Ker} f_i}$ it is not difficult to see that we get $f=\sum\limits_{i=1}^n f_i\cdot h_i$ so $f\in X$. %by the above considerations we can find This shows that for every $g_i\in {\rm Hom}^C(\frac{C}{{\rm Ker}f_i},C)$ there is $h_i\in R$ such that $g_i=\overline{f_i}\cdot h_i$ and therefore it follows that ${\rm Hom}^C(\frac{C}{{\rm Ker}f_i},C)$ is generated by $\overline{f_i}$ as a right $R$ module. By this and by the above exact sequence, we find that there are $g_i:C\rightarrow C$ such that $\varphi(\sum\limits_{i=1}^n\overline{f_i}\cdot g_i)=f$ and then by composing to the cannonocal projections we get that $f\in \sum\limits_{i=1}^{n}f_i\cdot R=X$. Therefor $(X^{\perp})^{\perp}=X$.
%\begin{diagram}
%0 & \rTo & \frac{C}{\bigcap\limits_{i=1}^{n}{\rm Ker}f_i} & \rTo & \bigoplus\limits_{i=1}^{n}\frac{C}{{\rm Ker}f_i} \\
%& & \dTo^{\overline{f}} & \ldDots & \\
%& & C &\\
%\end{diagram}
\end{proof}

\begin{proposition}\label{1.fin}
Let $C$ be a coalgebra such that ${Rat}(M)$ splits off in any finitely generated left $C^*$-module $M$. Then any indecomposable injective left $C$-comodule $E$ contains only finite dimensional proper subcomodules. 
\end{proposition}
\begin{proof}
Let $T$ be the socle of $E$; then $T$ is simple and $E=E(T)$ is the injective envelope of $T$. We show that if $K\subseteq E(T)$ is an infinite dimensional subcomodule then $K=E(T)$. Suppose $K\subsetneq E(T)$. Then there is a left $C$-subcomodule (right $C^*$-submodule) $K\subsetneq L \subset E(T)$ such that $L/K$ is finite dimensional. We have an exact sequence of left $C^*$-modules:
$$0\rightarrow (L/K)^* \rightarrow L^* \rightarrow K^* \rightarrow 0$$
As $L/K$ is a finite dimensional left $C$-comodule, we have that $(L/K)^*$ is a rational left $C^*$-module; thus ${Rat}(L^*)\neq 0$. Also $L^*$ is finitely generated as it is a quotient of $E(T)^*$ which is a direct summand of $C^*$. We have $L^*={Rat}(L^*)\oplus X$ for some left $C^*$-submodule $X$ of $L^*$. Then ${Rat}(L^*)$ is finitely generated because $L^*$ is, so it is finite dimensional. As $L$ is infinite dimensional by our assumption, we have $X\neq 0$. This shows that $L^*$ is decomposable and finitely generated, thus it has at least two maximal submodules, say $M,N$. We have an epimorphism $E(T)^*\stackrel{f}{\rightarrow} L^*\rightarrow 0$ and then $f^{-1}(M)$ and $f^{-1}(N)$ are distinct maximal $C^*$-submodules of $E(T)^*$. But by \cite{I}, Lemma 1.4, $E(T)^*$ has only one maximal $C^*$-submodule which is $T^\perp$, so we have obtained a contradiction. 
\end{proof}

Let $C_0$ be the coradical of $C$, the sum of all simple subcomodules of $C$. By \cite{DNR}, Section 3.1, $C_0$ semisimple coalgebra that is a direct sum of simple subcoalgebras $C_0=\bigoplus\limits_{i\in I}C_i$ and each simple subcoalgebra $C_i$ contains only one type of simple left (or right) $C$-comodule; moreover, any simple left (or right) $C$-comodule is isomorphic to one contained in a $C_i$. 

\begin{proposition}\label{1.simples}
Let $C$ be a coalgebra such that the rational part of every finitely generated left $C^*$ module splits off. Then there is only a finite number of isomorphism types of simple left $C$-comodules, equivalently, $C_0$ is finite dimensional.
\end{proposition}
\begin{proof}
By the above considerations, if $S_i$ is a simple left $C$-subcomodule of $C_i$, we have that $(S_i)_{i\in I}$ forms a set of representatives for the isomorphism types of simple left $C$-comodules. Let $\Ss$ be a set of representatives for the simple right $C$-comodules. Let $E(C_i)$ be an injective envelope of the left $C$-comodule $C_i$ included in $C$; then as $C_0$ is essential in $C$ we have $C=\bigoplus\limits_{i\in I}E(C_i)$ as left $C$-comodules or right $C^*$-modules. Then $C^*=\prod\limits_{i\in I}E(C_i)^*$ as left $C^*$-modules. As $S_i\subseteq E(C_i)$, we have epimorphisms of left $C^*$-modules $E(C_i)^*\rightarrow S_i^*\rightarrow 0$ and therefore we have an epimorphism of left $C^*$-modules $C^*\rightarrow \prod\limits_{i\in I}S_i^*\rightarrow 0$. But there is a one-to-one correspondence between left and right simple $C$-comodules given by $\{S_i\mid i\in I\}\ni S\mapsto S^*\in \Ss$. Hence there is an epimorphism $C^*\rightarrow\prod\limits_{S\in \Ss}S\rightarrow 0$, which shows that the left $C^*$-module $P=\prod\limits_{S\in \Ss}S$ is finitely generated (actually generated by a single element). But then as ${Rat}({}_{C^*}P)$ is a direct summand in $P$, we must have that ${Rat}({}_{C^*}P)$ is finitely generated, so it is finite dimensional. Therefore, as $\Sigma=\bigoplus\limits_{S\in\Ss}S$ is a rational left $C^*$-module which is naturally included in $P$, we have $\Sigma\subseteq {Rat}(P)$. This shows that $\bigoplus\limits_{S\in\Ss}S$ is finite dimensional so $I$ must be finite. This is equivalent to the fact that $C_0$ is finite dimensional, because each $C_i$ is a simple coalgebra, thus a finite dimensional one.
\end{proof}

We shall say that a coalgebra is left (right) finite Rat-splitting if the rational part of any finitely generated left (right) $C^*$-module splits off.

\begin{proposition}\label{1.prop}
Let $C$ be a left finite Rat-splitting coalgebra. Then the following assertions hold:\\
%\begin{itemize}
%\item[(i)] $C$ is artinian as left $C$ comodule (equivalently, as right $C^*$ module).
%\item[(ii)] $C^*$ is Noetherian as left $C^*$ module.
%\item[(iii)] $C$ has at most countable dimension.
%\item[(iv)] $C$ is injective as left $C^*$ module.
%\end{itemize}
(i) $C$ is artinian as left $C$-comodule (equivalently, as right $C^*$-module).\\
(ii) $C^*$ is left Noetherian.\\
(iii) $C$ has at most countable dimension.\\
(iv) $C$ is injective as left $C^*$-module.
\end{proposition}
\begin{proof}
(i) We have a direct sum decomposition $C=\bigoplus\limits_{i\in F}E(S_i)$ where $C_0=\bigoplus\limits_{i\in F}S_i$ is the decomposition of $C_0$ into simple left $C$-comodules and $E(S_i)$ are injective envelopes of $S_i$ contained in $C$. Then $F$ is finite as $C_0$ is finite dimensional. Also, by Proposition \ref{1.fin} the $E(S_i)$ 's are artinian as they contain only finite dimensional proper subcomodules, thus $C$ is an artinian left $C$-comodule. \\
(ii) Take $I$ a left ideal of $C^*$ and suppose it is not finitely generated; then we can find a sequence $(x_k)_{k}$ of elements of $I$ such that denoting $I_k={}_{C^*}<x_1,x_2,\dots,x_k>$, $x_{k+1}\notin I_k$. Then, corresponding to the ascending chain of left $C^*$ submodules of $C^*$, $I_1\subsetneq I_2\subsetneq \dots\subsetneq I_k\subsetneq \dots$ we have a descending chain of right $C^*$ submodules of $C$, $I_1^\perp\supseteq I_2^\perp\supseteq\dots\supseteq I_k^\perp\supseteq\dots$, which must be stationary as $C_{C^*}$ is artinian, so $I_k^\perp=I_{k+1}^\perp=\dots$. Then $(I_k^\perp)^\perp=(I_{k+1}^\perp)^\perp=\dots$ and then by Lemma \ref{1.perp} we get that $I_k=I_{k+1}$, and then $x_{k+1}\in I_{k+1}=I_k$ which is a contradiction. \\
(iii) For any $i\in F$, if $E(S_i)$ is infinite dimensional, we may inductively build a sequence $(x_k)_k$ of elements of $E(S_i)$ such that $x_{k+1}\notin\sum\limits_{j=1}^{k} x_j\cdot C^*$ for any $k$ because $\sum\limits_{j=1}^{k}x_j\cdot C^*$ is always a finite dimensional comodule. Then $E(S_i)=\sum\limits_{k}x_k\cdot C^*$, because $\sum\limits_{k}x_k\cdot C^*$ is an infinite dimensional subcomodule of $E(S_i)$ and one can apply Lemma \ref{1.fin}. Now, as each $x_k \cdot C^*$ is finite dimensional, the conclusion follows. \\
%(iv) If $E$ is the injective hull of $C$ as a left $C^*$-module; then $Rat(E)$ is a direct summand of $E$; next $C$ is injective as a comodule so it splits off in $Rat(E)$. Consequently $C$ is a direct summand of the injective $C^*$-module $E$, hence it is itself injective.
(iv) As $C$ is a finite coproduct of $E(S_i)$'s it is enough to prove that each $E(S_i)$ is injective and by \cite{I0} Lemma 2, it is enough to prove that $E=E(S_i)$ splits off in any left $C^*$ module $M$ such as $M/E$ is 1-generated, that is, it is generated by an element $\hat{x}\in M/E$. Let $H={Rat}(C^*\cdot x)$; then there is $X<C^*\cdot x$ such that $H\oplus X=C^*\cdot x$. Then $E+H$ is a rational $C^*$ module so $(E+H)\cap X=0$; also $M=C^*\cdot x+E$ so $(E+H)+X=M$, showing that $E+H$ is a direct summand in $M$. But as $E$ is an injective comodule, we have that $E$ splits off in $E+H$, thus $E$ must split in $M$ and the proof is finished.
\end{proof}

\section{Chain Coalgebras}
Let $C$ be a coalgebra and denote by $C_0\subseteq C_1\subseteq C_2\subseteq\dots$ the coradical filtration of $C$, that is, $C_0$ is the coradical of $C$, and $C_{n+1}\subseteq C$ such that $C_{n+1}/C_n$ is the socle of the right (or left) $C$-comodule $C/C_n$ for all $n\in\NN$. Then $C_n$ is a subcoalgebra of $C$ for all $n$, and the same $C_n$ is obtained whether we take the socle of the left $C$-comodule $C/C_n$ or of the right $C$-comodule $C/C_n$. Put $C_{-1}=0$ and $R=C^*$. By \cite{DNR} we have $\bigcup\limits_{n\in\NN}C_n=C$. 

\begin{definition}
We say that a coalgebra $C$ is a left (right) chain coalgebra if and only if the lattice of the left (right) subcomodules of $C$ is a chain, that is, any two subcomodules $M,N$ of $C$ are comparable (given two subsets $A,B$ of a set $X$ we say that $A$ and $B$ are comparable if either $A\subseteq B$ or $B\subseteq A$).
\end{definition}

The following result shows that this definition is left-right symmetric and also characterizes all chain coalgebras.

\begin{proposition}\label{chainco}
The following assertions are equivalent for a coalgebra $C$:\\
(i) $C$ is a right chain coalgebra.\\
(ii) $C_{n+1}/C_n$ is either $0$ or a simple (right) comodule for all $n\geq -1$.\\
(iii) $C$ is a left chain coalgebra.
In this case $C$ and $C_n$, $n\geq -1$ are the only subcomodules (left, right, two-sided) of $C$ and $C_n$ is finite dimensional for all $n$.
\end{proposition}
\begin{proof}
(ii)$\Rightarrow$(i) We prove that any subcomodule of $C$ must be equal either to one of the $C_n$'s or to $C$. Let $M$ be a right subcomodule of $C$ and suppose $M\neq C$ and $M\neq 0$. Then there is $n\geq 0$ such that $C_n\nsubseteq M$ and let $n$ the minimal natural number with this property. Then we must have $C_{n-1}\subseteq M$ by the minimality of $M$ and we show that $C_{n-1}=M$. Indeed, if $C_{n-1}\subsetneq M$ we can find a simple subcomodule of $M/C_{n-1}$. But then $C_{n-1}\neq C$ so $C_{n-1}\neq C_n$ and as $C_n/C_{n-1}$ is the only simple subcomodule of $C/C_n$ we find $C_n/C_{n-1}\subseteq M/C_{n-1}$, that is $C_n\subseteq M$, a contradiction. This also proves the last statement of the proposition.\\
%We have that $C_0=C_0/C_{-1}$ is simple, so there is only one simple right $C$ comodule; otherwise $C_0=0$ so $C=0$ and we are done. We prove that any subcomodule of $C$ must be equal either to one of the $C_n$'s or to $C$. Let $M$ be a right subcomodule of $C$. If $M=0$ then we are done; otherwise $M$ must contain a simple subcomodule of $C$ and as $s(C)=C_0$ we have that $C_0\subseteq M$. If $M=C_0$ then we are done. Otherwise $M/C_0$ is a nonzero subcomodule of $C/C_0$ so it must contain a simple subcomodule of $C/C_0$; but $C_1/C_0$ is the only simple subcomodule of $C/C_0$ so we get $C_1/C_0\subseteq M/C_0$, equivalently $C_1\subseteq M$. Again we may have $M=C_1$, or $C_1\subsetneq M$ and we continue this process. In the general case, assume we have found that $C_n\subseteq M$. If $C_n=M$ then the proof is finished; otherwise as $M/C_n$ is a nonzero subcomodule of $C/C_n$ it must contain a simple subcomodule of $C/C_n$ and as $C_{n+1}/C_n$ is the only simple right subcomodule of the nonzero comodule $C/C_n$ (nonzero because $M\neq C_n$) we again find that $C_{n+1}\subseteq M$. Obviously this process either stops at one point $n$ so $M=C_n$ or it never ends and we have $C_n\subset M$ for all $n$ and then $C=\bigcup\limits_{n\in\NN}C_n\subseteq M\subseteq C$ so $M=C$. \\
(i)$\Rightarrow$(ii) If $C_{n+1}/C_n$ is nonzero and it is not simple then we can find $S_1$ and $S_2$ two distinct simple modules contained in $C/C_n$. Then $S_1=M_1/C_n$, $S_2=M_2/C_n$ and $M_1\cap M_2=C_n$, $M_1\neq C_n$, $M_2\neq C_n$ because $S_1\cap S_2=0$ and $S_1$ and $S_2$ are distinct simple subcomodules of $C/C_n$. But this shows that neither $M_1\subsetneq M_2$ nor $M_2\subsetneq M_1$ which is a contradiction.\\
(ii)$\Leftrightarrow$(iii) is proved similarly.
\end{proof}

Denote $J=C_0^\perp$; by \cite{DNR} Lemma 2.5.7 and Corollary 3.1.10 we have that $J=J(C^*)$ (the Jacobson radical of $C^*$) and $(J^n)^\perp=C_{n-1}$, so $J^n\subseteq ((J^n)^{\perp})^\perp=C_{n-1}^\perp$. As $\bigcup\limits_{n\in\NN}C_n=C$, we see that $\bigcap\limits_{n\in\NN}J^n=0$. %Also for a left $C^*$-module $M$ denote by $J(M)$ its Jacobson radical.

\begin{definition}
We say that a coalgebra $C$ is left almost finite if the left regular comodule ${}^{C}C$ has only finite dimensional proper subcomodules. 
\end{definition}

By Proposition \ref{chainco} a chain coalgebra is left and right almost finite.

\begin{proposition}\label{2.2}
Let $C$ be an left almost finite coalgebra. Then $C^*$ is left Noetherian; moreover all nonzero left ideals of $C^*$ have finite codimension.
\end{proposition}
\begin{proof}
%We have that $C=E(C_0)$ as left comodules, thus by Proposition \ref{1.fin} we have that every left subcomodule of $C$ is finite dimensional (all $C_n$ are finite dimensional).  by Proposition \ref{1.prop} $I$ is finitely generated and 
Then if $I$ is a nonzero left ideal of $C^*$ take $f\neq 0$, $f\in I$. Then $F:=(C^*\cdot f)^\perp$ is a left coideal of $C$ so $F$ is finite dimensional. Then $F^\perp$ has finite codimension; but $C^*\cdot f=F^\perp$ by Proposition \ref{1.perp}. This shows that $I$ has finite codimension. Consequently, $C^*$ is left Noetherian. %Moreover, $C^*\cdot f$ has finite codimension in $I$ too. Then if $X$ is a complement of $C^*\cdot f$ in $I$ and $f_1,\dots,f_k$ is a $K$-basis in $X$ we get that $I$ is generated by $f,f_1,\dots,f_k$ as a left $C^*$-module. % $I^\perp=f\cdot C^*$ is a left $C^*$ submodule of $C$, thus a right subcomodule of $C$.   and by Proposition\ref{1.perp} $F=I^\perp\neq C$, so it must have finite dimension. Therefor $I=(I^\perp)^\perp$ has finite codimension in $C^*$.\\
\end{proof}
%We will also denote by $R=C^*$.
\begin{proposition}\label{2.5}
If $C$ is a chain coalgebra, then $J^n=C_{n-1}^\perp$ for all $n$ and $0$ and $J^n$, $n\geq 0$ are the only ideals (left, right, two-sided) of $C^*$. Consequently, $C^*$ is a chain algebra. %If $C$ is infinite dimensional, then $C^*$ is a domain.
\end{proposition}
\begin{proof}
If $I$ is a left ideal of $C^*$, then by Lemma \ref{1.perp} and Proposition \ref{2.2} we have $(I^\perp)^\perp=I$. By Proposition \ref{chainco}, $I^\perp=C$ or $I^\perp=C_{n-1}$ for some $n\geq 1$ and therefore $I=(I^\perp)^\perp=C_{n-1}^\perp$ or $I=0$. We prove by induction on $n$ that $J^n=C_{n-1}^\perp$, for all $n\geq 1$. For $n=1$, $J=C_0^\perp$ and assume $J^n=C_{n-1}^\perp$ for some $n\geq 2$. We again have  $J^{n+1}=((J^{n+1})^\perp)^{\perp}$ and as $(J^{n+1})^\perp=C_n$, we get $J^{n+1}=C_n^\perp$ and the proof is finished.
\end{proof}

For a left $C^*$-module $M$ denote by $T(M)$ the set of all torsion elements of $M$, that is, $T(M)=\{x\in M\mid ann_{C^*}x\neq 0\}$.
If $C$ is a finite dimensional coalgebra, then obviously any right $C$-comodule is rational and the category of right comodules coincides to that of the left $C^*$-modules. Then it is interesting to investigate the infinite dimensional case. We first consider a special kind of coalgebra:

\begin{proposition}\label{2.cor}
Let $C$ be an infinite dimensional left almost finite coalgebra and let $R=C^*$. Then for any left $R$ module $M$ we have ${Rat}(M)=T(M)$; moreover, $x\in{Rat}(M)$ if and only if $R\cdot x$ is finite dimensional. %Consequently, the splitting with respect to the comodules subcategory of ${}_{R}{\Mm}$ coincides with the splitting with respect to the clasical torsion.
\end{proposition}
\begin{proof}
If $x\in {Rat}(M)$ then $R\cdot x$ is finite dimensional and then ${\rm ann}_R x$ must be of finite codimension, thus nonzero as $R$ is infinite dimensional. Conversely, if $x\in T(M)$ and $x\neq 0$ then $I={\rm ann}_R x$ is a nonzero left ideal of $R$ so it must have finite codimension by Proposition \ref{2.2}. Then as $R/I\simeq R\cdot x$ we get that $R\cdot x$ is finite dimensional. Also $I^\perp\subset C$ is a finite dimensional left subcomodule of $C$ and thus the subcoalgebra $C^\prime$ of $C$ generated by $I^\perp$ is finite dimensional. Taking Proposition \ref{1.perp} into account, $I=(I^\perp)^\perp\supseteq C^\prime{}^\perp$. Then $C^\prime{}^\perp\cdot x=0$, and $R\cdot x$ becomes a left $R/C^\prime{}^\perp$-module. Now note that as $C^\prime$ is a finite dimensional coalgebra with $C^\prime{}^*\simeq C^*/C^\prime{}^\perp$, $R\cdot x$ has a structure of right $C^\prime$-module. So we have a map $\rho:R\cdot x\rightarrow R\cdot x\otimes C^\prime$, $\rho(c)=c_0\otimes c_1$ such that $h\cdot c=h(c_1)c_0$ for $h\in C^\prime{}^*$. But then if $\pi:C^*\rightarrow C^\prime{}^*$ is the cannonical projection $\pi(r)=r|_{C^\prime}$, for $r\in R$ and $c\in R\cdot x$ we have $r\cdot c=\pi(r)\cdot c=\pi(r)(c_1)c_0=r(c_1)c_0$, so we may regard $\rho$ as a $C$ comultiplication of $R\cdot x$, thus $x\in{Rat}(M)$.
\end{proof}

%We call a coalgebra $C$ colocal if its coradical $C_0$ is a simple left (and right) comodule. \\

We show that a chain coalgebra is a (left and right) finite splitting coalgebra. The proof of this can be done by a standard extension to the noncommutative case of the proof of the theorem on the structure of finitely generated modules over a PID, and obtain as a consequence the fact that $T(M)$ is a direct summand in $M$ for every finitely generated module $M$. However we can do this in a more direct way.

\begin{theorem}\label{2.Th}
If $C$ is a chain coalgebra, then $C$ a left and right finite splitting coalgebra.
\end{theorem}
\begin{proof}
%It is easy to see that the only subcomodules of $C$ are the $C_n$ and the only left (right) ideals of $R$ are $J^n$. Then the proof can be done by a standard extension to the noncommutative case of the proof of the theorem on the structure of finitely generated modules over a PID (principal ideal domain). We outline the proof here:\\
%We include here also a different proof:\\
First notice that every torsion-free $R$ finitely generated module $M$ is free: indeed if $x_1,\dots,x_n$ is a minimal system of generators, then if $\lambda_1x_1+\dots+\lambda_nx_n=0$ with $\lambda_i$ not all zero, we may assume that $\lambda_1\neq 0$. Without loss of generality we may also assume that $\lambda_1R\supseteq \lambda_iR,\,\forall i$ as any two ideals of $R$ are comparable by Proposition \ref{2.5}. Therefore we have $\lambda_i=\lambda_1s_i$ for some $s_i\in R$. Then $\lambda_1x_1+\lambda_1s_2x_2+\dots+\lambda_1s_nx_n=0$ implies $x_1+s_2x_2+\dots+s_nx_n=0$ as $M$ is torsionfree and $\lambda_1\neq 0$. Hence $x_1\in{}_{R}<x_2,\dots,x_n>$, contradicting the minimality of $n$. \\
Now if $M$ is any left $R$ module and $T=T(M)={Rat}(M)$ then $T(M/T(M))=0$. Indeed take $\hat{x}\in T(M/T(M))$ and put $I=ann_{C^*}\hat{x}\neq 0$ so $I$ has finite codimension and $I$ is a two-sided ideal by Proposition 2.5. By Proposition 2.2 $I$ is generated by some $h_1,\dots,h_k\in I$. Then if $y\in Ix$ we have $y=f\cdot x$, $f\in I$ so $f=\sum\limits_{i=1}^nr_i\cdot h_i$ and $y=f\cdot x=\sum\limits_{i=1}^nr_i\cdot h_ix$. Therefore $Ix$ is generated by $h_1x,\dots,h_nx$. Because $I=ann_{R}\hat{x}$ we have $Ix\subseteq T={Rat}(M)$ (we use Proposition \ref{2.cor}) and as $Ix$ is finitely generated rational we get that $Ix$ has finite dimension. We obviously have an epimorphism $\frac{R}{I}\rightarrow \frac{Rx}{Ix}$ which shows that $Rx/Ix$ is finite dimensional because $I$ has finite codimension in $R$. Therefore we get that ${\rm dim}(Rx)={\rm dim}(Rx/Ix)+{\rm dim}(Ix)<\infty$ so then by Proposition \ref{2.cor} we have that $Rx$ is rational, thus $x\in T$ so $\hat{x}=0$.\\
Now as $M/T$ is torsion-free, there are $x_1,\dots,x_n\in M$ whose images $\hat{x_1},\dots,\hat{x_n}$ in $M/T$ form a basis. Then it is easy to see that $x_1,\dots,x_n$ are linearly independent in $M$. Then if $X=Rx_1+\dots+Rx_n$ we have $X+T=M$ and $X\cap T=0$, because if $a_1x_1+\dots+a_nx_n\in T\}$ we get $a_1\hat{x_1}+\dots+a_n\hat{x_n}=\hat{0}$ so $a_i=0,\,\forall i$ because $\hat{x_1},\dots,\hat{x_n}$ are independent in $M/T$. Thus $T(M)$ splits off in $M$ and the theorem is proved, as $T(M)={Rat}_R(M)$ by \ref{2.cor}.
%THE PROOF IS A STANDARD EXTENSION TO THE NONCOMMUTATIVE CASE OF THE PROOF OF THE THEOREM ON THE STRUCTURE OF FINITELY GENERATED MODULES OVER A PID (PRINCIPAL IDEAL DOMAIN); THE DETAILS TO BE COMPLETED.
\end{proof}

We will denote by $K_n$ the coalgebra with a basis $c_0,c_1,\dots,c_{n-1}$ and comultiplication $c_k\mapsto \sum\limits_{i+j=k}c_i\otimes c_j$ and counit $\varepsilon(c_i)=\delta_{0,i}$. The coalgebra $\bigcup\limits_{n\in \NN}K_n$ having a basis $c_n,n\in \NN$ and comultiplication and counit given by these equations is called the divided power coalgebra (see \cite{DNR}).% and will be denoted by $K^\prime$.

\begin{lemma}\label{2.7}
Let $C$ be a finite dimensional chain coalgebra over an algebraically closed field. Then $C$ is isomorphic to $K_n$ for some $n\in\NN$. 
\end{lemma}
\begin{proof}
Let $A=C^*$; we have ${\rm dim}\,C_0=1$ because $K$ is algebraically closed (thus ${\rm End}_AC_0$ is a skewfield containing $K$). Thus ${\rm dim}\,C_k=k$ for all $k$ for which $C_k\neq C$. As $C^*$ is finite dimensional $J^n=0$ for some $n$ and let $n$ be minimal with this property. By Proposition \ref{2.5} $J^k=C_{k-1}^\perp$. Then $J^k/J^{k+1}$ has dimension equal to the dimension of $C_{k}/C_{k-1}$ which is $1$ for $k<n$, because $C_{k+1}/C_k$ it is a simple comodule isomorphic to $C_0$. We then have that $J^k/J^{k+1}$ is generated by any of its nonzero elements. Choose $x\in J\setminus J^2$. We prove that $x^{n-1}\neq 0$. Suppose the contrary holds and take $y_1,\dots,y_{n-1}\in J$. As $x$ generates $J/J^2$, there is $\lambda\in K$ such that $y_{1}-\lambda x\in J^2$ and then $y_{1}x^{n-2}-\lambda x^{n-1}\in J^n$, so $y_{1}x^{n-2}\in J^n=0$ because $x^{n-1}=0$. Again, there is $\mu\in K$ such that $y_{2}-\mu x\in J^2$ and then $y_{1}y_2-\mu y_1 x\in J^3$ so $y_1y_2x^{n-3}\in J^n$ ($y_1x^{n-2}=0$). By continuing this procedure, one gets that $y_1y_2\dots y_{n-2}x=0$ and then we again find $\alpha\in K$ with $y_{n-1}-\alpha x\in J^2$, thus $y_1\dots y_{n-1}-\alpha y_1\dots y_{n-2}x\in J^n=0$. This shows that $y_1\dots y_{n-1}=0$ for all $_1\dots y_{n-1}=0$. Thus $J^{n-1}=0$, a contradiction. \\
As $x^{n-1}\neq 0$ we see that $x^k\in J^k\setminus J^{k+1}$ for all $k=0,\dots,n-1$, so $J^k/J^{k+1}$ is generated by the class of $x^k$. Now if $y\in A$, there is $\lambda_0\in K$ such that $y-\lambda_0\cdot 1_A\in J$ (either $y\in J$ or $y$ generates $A/J$). As $J/J^2$ is 1 dimensional and generated by the image of $x$, there is $\lambda_1\in K$ such that $y-\lambda_0-\lambda_1x\in J^2$. Again, as $J^2/J^3$ is 1 dimensional generated by the image of $x^2$, there is $\lambda_2\in K$ such that $y-\lambda_0-\lambda_1x-\lambda_2x^2\in J^3$. By continuing this procedure we find $\lambda_0,\dots,\lambda_{n-1}\in K$ such that $y-\lambda_0-\lambda_1x-\dots-\lambda_{n-1}x^{n-1}\in J^{n}=0$, so $y=\lambda_0+\lambda_1x+\dots+\lambda_{n-1}x^{n-1}$. This obviously gives an isomorphism between $A$ and $K[X]/(X^n)$. Therefore $C$ is isomorphic to $K_n$, because there is an isomorphism of $K$-algebras $K_n^*\simeq K[X]/(X^n)$.%sheds light on an isomorphism between $A$ and $K[X]/(X^n)$, thus $C$ is isomorphic to $K_n$, the coalgebra dual to the algebra $K[X]/(X^n)$. 
\end{proof}

\begin{theorem}\label{2.ThK}
If $K$ is an algebraically closed field and $C$ is an infinite dimensional chain coalgebra, then $C$ is isomorphic to the divided power coalgebra.
\end{theorem}
\begin{proof}
By the previous Lemma we have that $C_n\simeq K_n$ for all $n$. If $e\in C_0$, $\Delta(e)=\lambda e\otimes e,\,\lambda\in K$, then for $c_0=\lambda e$ we get $\Delta(c_0)=c_0\otimes c_0$. Suppose we constructed a basis $c_0,c_1,\dots,c_{n-1}$ for $C_{n-1}$ with $\Delta(c_k)=\sum\limits_{i+j=k}c_i\otimes c_j$, $\varepsilon(c_i)=\delta_{0,i}$. Denote by $A_n=C_n^*$ the dual of $C_n$; for the rest of this proof, if $V\subseteq C_n$ is a subspace of $C_n$ we write $V^\perp$ for the set of the functions of $A_n$ which are 0 on $V$. %We choose $d_n\in C_{n}\cap{\rm Ker}\varepsilon\setminus C_{n-1}$ and denote by $A_n=C_n^*$ the dual of $C_n$; pick $E_0,E_1,\dots,E_{n-1},D_n$ a dual basis in $A_n$ for $c_0,c_1,\dots,c_{n-1},d_n$. let $E_1^n(d_n)=\lambda_n$ and $c_n=\frac{1}{\lambda_n}d_n$ It is easy to see that 
Choose $E_1\in C_0^\perp\setminus C_1^\perp$; then $E_1^n\neq 0$ and $E_1^{n+1}=0$ as in the proof of Lemma \ref{2.7} ($E_1\in A_n$). This shows that $E_1^k\in C_{k-1}^\perp\setminus C_k^\perp$, that $\varepsilon\vert_{C_n},E_1,\dots,E_1^n$ exhibits a basis for $A_n$ and that  there is an isomorphism of algebras $A_n\simeq K[X]/(X^{n+1})$ taking $E_1$ to $\hat{X}$. We can easily see that $E_1^i(c_j)=\delta_{ij}$, $\forall k=0,1,\dots,n-1$ and then by a standard linear algebra result we can find $c_n\in C_n$ such that $E_1^n(c_n)=1$ and $E_1^n(c_i)=0$ for $i<n$. Then by dualization, the relations $E_1^i(c_j)=\delta_{ij}$, $\forall i,j=0,1,\dots,n$ become $\Delta(c_k)=\sum\limits_{i+j=k}c_i\otimes c_j$, $\forall k=0,1,\dots,n$. Therefore we may inductively build the basis $(c_n)_{n\in\NN}$ with $\varepsilon(c_k)=\delta_{0k}$ and $\Delta(c_n)=\sum\limits_{i+j=n}c_i\otimes c_j,\,\forall n$.
\end{proof}

In the following we construct an example of a chain coalgebra that is not cocommutative and thus different of the divided power coalgebra over $K$. Recall that if $A$ is a $k$ algebra, $\varphi:A\rightarrow A$ is a morphism and $\delta:A\rightarrow A$ is a $\varphi$-derivation (that is a linear map such that $\delta(ab)=\delta(a)b+\varphi(a)\delta(b)$ for all $a,b\in A$), we may consider the Ore extension $A[X,\varphi,\delta]$ which is $A[X]$ as a vector space and with multiplication induced by $Xa=\varphi(a)X+\delta(a)$.
%\begin{example}
Let $K$ be a subfield of $\mathbb R$, the field of real numbers. Let $D$ be the subalgebra of Hamilton's quaternion algebra having the set $B=\{1,\overline{i},\overline{j},\overline{k}\}$ as a vector space basis over $K$. Recall that multiplication is given by the rules $\overline{i}\cdot\overline{j}=-\overline{j}\cdot\overline{i}=\overline{k}$; $\overline{j}\cdot\overline{k}=-\overline{k}\cdot\overline{j}=\overline{i}$; $\overline{k}\cdot\overline{i}=-\overline{i}\cdot\overline{k}=\overline{j}$; $\overline{i}^2=\overline{j}^2=\overline{k}^2=-1$. Denote by $\sigma:D\rightarrow D$ the linear map defined on the basis of $D$ by
$$
\sigma=\left(
\begin{array}{cccc}
1 & \overline{i} & \overline{j} & \overline{k}\\
1 & \overline{j} & \overline{k} & \overline{i}
\end{array}
\right)
$$
It is not difficult to see then that $\sigma$ is an algebra automorphism, and that $D$ is a division algebra (skewfield). Our example will be such an Ore extension constructed with a trivial derivation: denote by $D_{\sigma}[X]=D[X,\sigma,0]$ the Ore extension of $D$ constructed by $\sigma$ with the derivation $\varphi$ equal to $0$ everywhere. Then a basis for $D_{\sigma}[X]$ over $K$ consists of the elements $uX^k$, with $u\in B$ and $k\in\NN$. Also denote by $A_n=D_{\sigma}[X]/<X^n>$ the algebra obtained by factoring out the two-sided ideal generated by $X^n$ from $D_{\sigma}[X]$.
%\end{example}
\begin{proposition}\label{An}
The two sided ideal $<X^n>$ of $D_\sigma[X]$ consists of elements of the form $f=\sum\limits_{l=n}^{n+m}a_lX^l$. Moreover, the only (left, right, two-sided) ideals containing $<X^n>$ are the ideals $<X^l>$, $l=0,\dots,n$ and consequently $A_n$ is a chain $K$ algebra.
\end{proposition}
\begin{proof}
It is clear by the multiplication rule $Xa=\sigma(a)X$ for $a\in B$ that elemens of $D_\sigma[X]$ are of the type $\sum\limits_{l=0}^Na_lX^l$ and that every element of $A_n$ is a "polynomial" of the form $f=a_0+a_1x+\dots+a_{n-1}x^{n-1}$, with $a_l\in D$ and where $x$ represents the class of $X$. Such an element $f$ is invertible if and only if $a_0\neq 0$. To see this, first note that if $a_0=0$ then $f$ is nilpotent, as $x$ is nilpotent and one has $f^l\in<x^l>$ by successively using the relation $xa=\sigma(a)x$. Conversely write $f=a_0\cdot(1+a_0^{-1}a_1x+\dots+a_0^{-1}a_{n-1}x^{n-1})$ and note that the element $g=a_0^{-1}a_1x+\dots+a_0^{-1}a_{n-1}x^{n-1}$ is nilpotent as before, so $1+g$ must be invertible in $A_n$ and therefore $f$ must be invertible. Thus we may write every element $f=a_lx^l+...a_{n-1}x^{n-1}$ of $A_n$ as the product $f=(a_l+a_{l+1}x+\dots+a_{n-1}x^{n-1-l})\cdot x^l=g\cdot x^l$ with invertible $g$. Then if $I$ is a left ideal of $A_n$ and $f\in I$, we have $f=g\cdot x^l$ for an invertible element $g$ and some $l\leq n$. Hence it follows that $x^l\in I$. Taking the smallest number $l$ with the property $x^l\in I$, we obviously have that $I=<x^l>$.
\end{proof}

Let $C_n$ denote the coalgebra dual to $A_n$. Note that $A_n$ has a $K$ basis $\Bb=\{ax^l\mid a\in B,\,l\in 0,1,\dots,n-1\}$ and we have the relations $(ax^i)(bx^j)=a\sigma^i(b)x^{i+j}$. Let $(E_i^a)_{a\in B,i\in\overline{0,n-1}}$ the basis of $C_n$ which is dual to $\Bb$, that is, $E_i^a(bx^j)=\delta_{ij}\delta_{ab}$ for all $a,b\in B$ and $i,j\in {\NN}$. Also, for $i\in\NN$ and $a\in B$ denote by $i\cdot a=\sigma^i(a)$, the action of $\NN$ on $B$ induced by $\sigma$.

\begin{proposition}
With the above notations, denoting by $\Delta_n$ and $\varepsilon_n$ the comultiplication and respectively, the counit of $C_n$ we have
\begin{eqnarray*}
\Delta_n(E_p^c)=\sum\limits_{i+j=p;\,\,a(i\cdot b)=\pm c}c^{-1}a(i\cdot b)E_i^a\otimes E_j^b
\end{eqnarray*}
and
\begin{eqnarray*}
\varepsilon_n(E_p^c)=\delta_{p,0}\delta_{c,1}.
\end{eqnarray*}
\end{proposition}
\begin{proof}
For $u,v\in B$ and $k,l\in\NN$ we have $E_p^c(ux^k\cdot vx^l)=E_p^c(u(k\cdot v)x^{k+l})$ and as $k\cdot v\in B$ by the formulas defining $D$ we have that if $d=u(k\cdot v)$ then either $d\in B$ or $-d\in B$. Then $E_p^c(ux^k\cdot vx^l)=E_p^c(dx^{k+l})=\delta_{k+l,p}\delta_{u(k\cdot v),\pm c}c^{-1}u(k\cdot v)$ as the sign of this expression must be $1$ if $d\in B$ and $-1$ if $d\notin B$, and this is exactly $c^{-1}u(k\cdot v)$ when $u(k\cdot v)=\pm c$. We also have
\begin{eqnarray*}
\sum\limits_{i+j=p;\,\,a(i\cdot b)=\pm c}c^{-1}a(i\cdot b)E_i^a(ux^k)E_j^b(vx^l) & = & \sum\limits_{i+j=p;\,\,a(i\cdot b)=\pm c}\delta_{k,i}\delta_{u,a}\delta_{l,j}\delta_{v,b}c^{-1}a(i\cdot b)\\
& = & \delta_{k+l,p}\delta_{u(k\cdot v),\pm c}c^{-1}u(k\cdot v)
\end{eqnarray*}
and therefore we get 
$$\sum\limits_{i+j=p;\,\,a(i\cdot b)=\pm c}c^{-1}a(i\cdot b)E_i^a(ux^k)E_j^b(vx^l)=E_p^c(ux^k\cdot vx^l) $$
As this is true for all $ux^k, vx^l\in\Bb$, by the definition of the comultiplication of the coalgebra dual to an algebra, we get the first equality in the statement of the proposition. The second one is obvious, as $\varepsilon_n(E_p^c)=E_p^c(1\cdot X^0)=\delta_{p,0}\delta_{c,1}$.
\end{proof}

Now notice that there is an injective map $C_n\subset C_{n+1}$ taking $E_i^c$ from $C_n$ to $E_i^c$ from $C_{n+1}$. Therefore we can regard $C_n$ as subcoalgebra of $C_{n+1}$. Denote by $C=\bigcup\limits_{n\in \NN}C_n$; it has a basis formed by the elements $E_n^c,\,n\in\NN,\,c\in B$ and comultiplication $\Delta$ and counit $\varepsilon$ given by 
$$\Delta(E_n^c)=\sum\limits_{i+j=n;\,\,a(i\cdot b)=\pm c}c^{-1}a(i\cdot b)E_i^a\otimes E_j^b$$
and
$$\varepsilon(E_n^c)=\delta_{n,0}\delta_{c,1}.$$
By Proposition \ref{An} we have that $A_n$ is a chain algebra and therefore $C_n=A_n^*$ is a chain coalgebra. Therefore, we get that the coradical filtration of $C$ is $C_0\subseteq C_1\subseteq C_2\subseteq\dots$ and that this is a chain coalgebra which is obviously non-cocommutative.

\section{The co-local case}

Throughout this section we will assume (unless otherwise specified) that $C$ is left finite Rat-splitting and that it is a colocal coalgebra, that is, $C_0$ is a simple left (and consequently simple right) $C^*$-module. Then as $J=C_0^\perp$, $C^*$ is a local algebra. We will also assume that $C$ is not finite dimensional, thus by Proposition \ref{1.prop} $C$ has a countable basis. %We have the coradical filtration of $C$, $C_0\subseteq C_1\subseteq\dots\subseteq C_n\subseteq\dots$. Denote $J=C_0^\perp$; by \cite{DNR} we have that $J=J(C^*)$ (the Jacobson radical of $C^*$) and $J^n=C_n^\perp$. As $\bigcup\limits_{n\in\NN}C_n=C$, se see that $\bigcap\limits_{n\in\NN}J^n=0$. 
We have that $C$ is the injective envelope of $C_0$ as left comodules, thus by Proposition \ref{1.fin} we have that every left subcomodule of $C$ is finite dimensional (all $C_n$ are finite dimensional). Then if $I$ is a left ideal of $C^*$ different from $C^*$, by Proposition \ref{1.prop} and Proposition \ref{2.2} $I$ is finitely generated and of finite codimension. %and $I^\perp$ is a right $C^*$ submodule of $C$, thus a left subcomodule of $C$,  and by Proposition\ref{1.perp} $I^\perp\neq C$, so it must have finite dimension. Therefor $I=(I^\perp)^\perp$ has finite codimension in $C^*$.\\
Denote again $R=C^*$. Also for a left $R$-module $M$ denote by $J(M)$ the Jacobson radical of $M$.

\begin{proposition}\label{2.domain}
With the above notations, $R$ is a domain.
\end{proposition}
\begin{proof}
Let $S={\rm End}({}^CC,{}^CC)$. Note that $S$ is a ring with multiplication equal to the composition of morphisms and that $S$ is isomorphic to $R$ by an isomorphism that takes every morphism of left $C$-comodules $f\in S$ to the element $\varepsilon\circ f\in R$. Then it is enough to show that $S$ is a domain. If $f:C\rightarrow C$ is a nonzero morphism of left $C$ comodules, then ${\rm Ker}(f)\subsetneq C$ is a proper left subcomodule of $C$ so it must be finite dimensional. Then as $C$ is not finite dimensional we see that ${\rm Im}(f)\simeq C/{\rm Ker}(f)$ is an infinite dimensional subcomodule of $C$. Thus ${\rm Im}(f)=C$, and therefore every nonzero morphism of left comodules from $C$ to $C$ must be surjective. Now if $f,g\in S$ are nonzero then they are surjective so $f\circ g$ is surjective and thus $f\circ g\neq 0$. 
\end{proof}

\begin{proposition}\label{2.accp}
$R$ satisfies ACCP on right ideals and also on left ideals.
\end{proposition}
\begin{proof}
Suppose there is an ascending chain of right ideals $x_0\cdot R\subsetneq x_1\cdot R\subsetneq x_2\cdot R\subsetneq\dots$ that is not stationary. Then there are $(\lambda_n)_{n\in\NN}$ in $R$ such that $x_n=x_{n+1}\cdot\lambda_{n+1}$. Note that $\lambda_{n+1}\in J$, because otherwise $\lambda_{n+1}$ would be invertible in $R$ as $R$ is local and then we would have $x_{n+1}=x_n\cdot \lambda_n^{-1}$. This would yield $x_n\cdot R=x_{n+1}\cdot R$, a contradiction. Then $x_1=x_{n+1}\cdot\lambda_{n+1}\lambda_{n}\dots\lambda_2$, so $x_1\in J^n$ for all $n\in\NN$, showing that $x_1\in\bigcap\limits_{n\in \NN}J^n=0$. Thus we obtain a contradiction: $x_0\cdot R\subsetneq x_1\cdot R=0$. The statement is obvious for left ideals as ${}_{R}R$ is Noetherian.
\end{proof}

%We will say that two subsets $A$ and $B$ of $R$ are comparable if $A\subseteq B$ or $B\subseteq A$. 
The next proposition contains the main idea of the result.

\begin{proposition}\label{2.key}
Suppose $\alpha R$ and $\beta R$ are two right ideals that are not comparable. Then any two principal right ideals of $R$ contained in $\alpha R\cap \beta R$ are comparable.
\end{proposition}
\begin{proof}
Take $aR,bR\subseteq \alpha R\cap \beta R$, so $a=\alpha x=\beta y$ and $b=\alpha u=\beta v$; we may obviously assume that $a,b\neq 0$ as otherwise the assertion is obvious. Then $\alpha,\beta,x,y,u,v$ are nonzero. Denote by $L$ the left submodule of $R\times R$ generated by $(x,u)$ and by $M$ the quotient module $\frac{R\times R}{L}$. We write $\overline{(s,t)}$ for the image of the element $(s,t)$ through the canonical projection $\pi:R\times R\rightarrow M$. We have $\overline{(y,v)}\neq\overline{(0,0)}$ as otherwise $(y,v)=\lambda(x,u)$ for some $\lambda\in R$; then we would have $y=\lambda x$, $v=\lambda u$ so $\beta y=\beta\lambda x=\alpha x$ and then $\beta\lambda=\alpha$ (because $R$ is a domain), a contradiction to $\alpha R\subsetneq \beta R$. Also $\beta\cdot\overline{(y,v)}=\alpha\cdot\overline{(x,u)}=\overline{(0,0)}$ with $\beta\neq 0$. This shows that $\overline{(0,0)}\neq\overline{(y,v)}\in T=T(M)$, so $T(M)\neq 0$. Take $X<M$ such that $M=T\oplus X$. We must have $X\neq 0$, as otherwise $\overline{(1,0)}\in T$ so there would be a nonzero $\lambda\in R$ and a $\mu\in R$ such that $\lambda\cdot(1,0)=\mu\cdot(x,u)\in L$. But then $\lambda=\mu x$, $0=\mu u$, so $\mu=0$ ($u\neq 0$) showing that $\lambda=0$, a contradiction. \\
Now note that $x$ and $u$ are not invertible, as otherwise, for $x$ invertible, $\alpha x=\beta y$ implies $\alpha\in\beta R$ so $\alpha R\subseteq \beta R$; the same can be inferred if $u$ is invertible. Therefore $x,u\in J$ as $R$ is local so $L\subseteq J\times J$. Hence $J(M)=J\times J/L$ so $M/J(M)=\frac{R\times R/L}{J\times J/L}\simeq R\times R/J\times J$ which has dimension 2 as a module over the skewfield $R/J$. Since $M=T\oplus X$ and $M$ is finitely generated, then so are $T$ and $X$ and therefore $J(X)\neq X$ and $J(T)\neq T$. Then as $\frac{M}{J(M)}=\frac{T}{J(T)}\oplus\frac{X}{J(X)}$ has dimension 2 over $R/J$, it follows that both $T/J(T)$ and $X/J(X)$ are simple. Hence $T$ and $X$ are local, and as they are finitely generated, it follows that they are generated by any element not belonging to their Jacobson radical. Let $T'$ (respectively $X'$) be the inverse images of $T$ (and $X$ respectively) in $R\times R$ and $t\in T'$ and $s\in X'$ be such that $Rt+L=T'$ and $Rs+L=X'$. We have $R\times R=T'+X'=Rt+L+Rs+L=(Rt+Rs)+L\subseteq (Rt+Rs)+J\times J\subseteq R\times R$ so $(Rt+Rs)+J\times J=R\times R$. Therefore we obtain $Rt+Rs=R\times R$ because $J\times J$ is small in $R\times R$.\\
%As $X$ and $T$ are finitely generated $R$ modules (as direct summands in $R\times R$) we have that $J(X)\neq X$ and $J(T)\neq T$. We show that they are local modules. Assume they are not; then there are proper submodules $T_1$, $X_1$ of $T$ and $X$ respectively such that the modules $X/X_1$ and $T/T_1$ are non-zero semisimple and at least one of them is not simple. Denoting by $U$ the inverse image of $X_1\oplus T_1$ through $\pi$, we have $$X/X_1\oplus T/T_1=(X\oplus T)/(X_1\oplus T_1)=\frac{(R\times R)/L}{U/L}\simeq R\times R/U$$ and then $R\times R/U$ must be semisimple of length at least $3$; but $J(R\times R)=J\times J\subseteq U$ so $J(R\times R)$ has co-length 2 and we obtain a contradiction. As $T$ is finitely generated and local, we see that $T$ must be generated by some $\overline{t}\in T$. Indeed, if $\overline{t}\in T\setminus J(T)$, then $T=R\overline{t}$, because otherwise as $T$ is finitely generated, $R\overline{t}$ is contained in a maximal submodule of $T$ so $R\overline{t}\subseteq J(T)$ because $J(T)$ is the only maximal $R$-submodule of $T$. Similarly, $X$ is generated by some $\overline{s}\in X$. Denote $T'$ (respectively $X'$) the inverse image of $T$ (respectively $X$) in $R\times R$.\\
%Now note that $x$ and $u$ are not invertible, as otherwise, for $x$ invertible, $\alpha x=\beta y$ implies $\alpha\in\beta R$ so $\alpha R\subseteq \beta R$; the same can be inferred if $u$ invertible. Therefore $x,u\in J$ as $R$ is local so $L\subseteq J\times J$. 
Write $t=(p,q)\in T'$. Then $\overline{t}=t+L\in T$ implies that there is $\lambda\neq 0$ in $R$ such that $\lambda\overline{t}=\overline{0}\in M$ and therefore there is $\mu\in R$ with $\lambda(p,q)=\mu(x,u)$. We show that either $p\notin J$ or $q\notin J$. Indeed assume otherwise: $t=(p,q)\in J\times J$. Then we get $Rt\subseteq J\times J$. Because $Rt+Rs=R\times R$ we see that $R\times R/J\times J$ must be generated over $R$ by the image of $s$. This shows that the $R/J$ module $R\times R/J\times J=(R/J)^2$ has dimension $1$ and this is obviously a contradiction.\\
Finally, suppose $p\notin J$ so $p$ is invertible; then the equations $\lambda p=\mu x$, $\lambda q=\mu u$ imply $\lambda=\mu xp^{-1}$ and $\mu xp^{-1}q=\mu u$. But $\mu\neq 0$ because $p$ is invertible and $\lambda\neq 0$. Therefore we obtain $u=xp^{-1}q$; thus $b=\alpha u=\alpha xp^{-1}q=ap^{-1}q$ showing that $b\in aR$ i.e. $bR\subseteq aR$. Similarly if $q$ is invertible, we get $aR\subseteq bR$.
\end{proof}

%The main result of the paper is the following:

\begin{theorem}\label{3.Th}
If $C$ is a left finite splitting (infinite dimensional) local coalgebra, then $C$ is a chain coalgebra.%every factor of the coradical filtration $C_{n+1}/C_n$ is simple.
\end{theorem}
\begin{proof}
We first show that every two principal left ideals of $R$ are comparable. Suppose there are two left ideals of $R$, $R\cdot x_0$ and $R\cdot y_0$ that are not comparable. Then as they have finite codimension and $C^*$ is infinite dimensional, we have $Rx_0\cap Ry_0\neq 0$ and take $0\neq \alpha x_0=\beta y_0\in Rx_0\cap Ry_0$. Then the right ideals $\alpha R$ and $\beta R$ are not comparable, as otherwise, if for example $\alpha R\subseteq \beta R$, we would have a relation $\alpha=\beta\lambda$ so $\alpha x_0=\beta\lambda x_0=\beta y_0$. As $\beta\neq 0$ we get $\lambda x_0=y_0$ because $R$ is a domain, and then $R y_0\subseteq R x_0$, a contradiction.\\
By Proposition \ref{2.accp} the set $\{\lambda R\mid \lambda R\subseteq \alpha R\cap\beta R\}$ is Noetherian (relative to inclusion) and let $\lambda R$ be a maximal element. If $x\in \alpha R\cap \beta R$ then by Proposition \ref{2.key} we have that $xR$ and $\lambda R$ are comparable and by the maximality of $\lambda R$ it follows that $xR\subseteq \lambda R$, so $x\in\lambda R$. Therefore $\alpha R\cap \beta R=\lambda R$. Note that $\lambda\neq 0$, because $\alpha R$ and $\beta R$ are nonzero ideals of finite codimension. Then we see that $\lambda R\simeq R$ as right $R$ modules, because $R$ is a domain, and again by Proposition \ref{2.key} any two principal right ideals of $\lambda R=\alpha R\cap \beta R$ are comparable, so the same must hold in $R_{R}$. But this is in contradiction with the fact that $\alpha R$ and $\beta R$ are not comparable, and therefore the initial assertion is proved. \\
Now we prove that $J^n/J^{n+1}$ is a simple right module for all $n$. As $R/J$ is semisimple (it is a skewfield) and $J^n/J^{n+1}$ has an $R/J$ module structure, it follows that $J^n/J^{n+1}$ is a semisimple left $R/J$-module and then $J^n/J^{n+1}$ is semisimple also as $R$-module. If we assume that it is not simple, then there are $f,g\in J^n\setminus J^{n+1}$ such that $R\hat{f}=(Rf+J^{n+1})/J^{n+1}$ and $R\hat{g}=(Rg+J^{n+1})/J^{n+1}$ are different simple $R$-modules, so $R\hat{f}\cap R\hat{g}=\hat{0}$ in $J^n/J^{n+1}$. Then $(Rf+J^{n+1})\cap(Rg+J^{n+1})=J^{n+1}$ which shows that $Rf$ and $Rg$ cannot be comparable, a contradiction. 
%We can see that because any two principal left submodules in $R$ are comparable, the same holds for the left $R$-module $J^n/J^{n+1}$ and as $J^n/J^{n+1}$ is semisimple, this can only happen if it is a simple $R$-module. Therefore $J^n/J^{n+1}\simeq C_0$ as left $R$-modules. 
As $J^n=C_{n-1}^\perp$, we see that ${\rm dim}(C_{n-1})={\rm codim}(J^n)$. Then for $n\geq 1$, ${\rm dim}(C_{n}/C_{n-1})={\rm dim}(C_n)-{\rm dim}(C_{n-1})={\rm codim}_R(J^n)-{\rm codim}_R(J^{n+1})={\rm dim}(J^n/J^{n+1})={\rm dim}(C_0)$. Because $C_0$ is the only type of simple right $C$-comodule, this last relation shows that the right $C$-comodule $C_n/C_{n-1}$ must be simple. Therefore $C$ must be a chain coalgebra. 
%We first prove that any two 
\end{proof}

We may now combine the results of Sections 2 and 3 and obtain

\begin{corollary}
Let $C$ be a co-local coalgebra. Then $C$ is a left (right) finite splitting coalgebra if and only if $C$ is a chain coalgebra. Moreover, if the base field $K$ is algebraically closed then this is further equivalent to the fact that $C$ is isomorphic to the divided power coalgebra.
\end{corollary}
\begin{proof}
This follows from Theorems \ref{2.Th}, \ref{2.ThK} and \ref{3.Th}.
\end{proof}

\bigskip\bigskip\bigskip

\begin{center}
\sc Acknowledgment
\end{center}
The author wishes to thank his Ph.D. adviser C. N\u ast\u asescu for very useful remarks on the subject as well as for his continuous support throughout the past years.

\vspace*{3mm} 
\begin{flushright}
\begin{minipage}{148mm}\sc\footnotesize

Miodrag Cristian Iovanov\\
University of Bucharest, Faculty of Mathematics, Str. Academiei 14\\ 
RO-010014, Bucharest, Romania\\
{\it E--mail address}: {\tt
yovanov@walla.com}\vspace*{3mm}

\end{minipage}
\end{flushright}

\newpage
\bigskip\bigskip\bigskip
\begin{thebibliography}{J\c{S}}

%\bibitem[Mo]{Mo}
%K. Morita, \emph{Adjoint pairs of functors and Frobenius extensions}, Sci. Rep. Tokyo Kyoiku Daigaku Sect. A, 9, 40-71 (1965).


%\bibitem[DNR]{DNR} 
%S. D\u asc\u alescu, C. N\u ast\u asescu, \c S. Raianu, \emph{Hopf Algebras: an introduction}. Vol. 235. Pure and Applied Mathematics, Marcel Dekker, New York, 2001.

\bibitem[AN]{AN}
T. Albu, C. N\u ast\u asescu, \emph{Relative Finiteness in Module Theory}, Monogr. Textbooks Pure Appl. Math., vol. 84, Dekker, New York 1984.

\bibitem[AF]{AF}
D. Anderson, K.Fuller, \emph{Rings and Categories of Modules}, Grad. Texts in Math., Springer, Berlin-Heidelberg-New York, 1974.

\bibitem[D1]{Doi2} 
Y. Doi, \emph{Homological Coalgebra}, J. Math. Soc. Japan \textbf{33}(1981), 31-50.

\bibitem[BW]{BW} 
T. Brzezi\'nski and R. Wisbauer, \emph{Corings and comodules}, London Math. Soc. Lect. Notes Ser. {\bf 309}, Cambridge University Press, Cambridge, 2003.

\bibitem[CIDN]{NT0}
F. Casta$\rm\tilde{n}$o Iglesias, S. D\u asc\u alescu, C. N\u ast\u asescu, \emph{Symmetric Coalgebras}, J. Algebra {\bf 279} (2004) 326-344.

\bibitem[DNR]{DNR} 
S. D\u asc\u alescu, C. N\u ast\u asescu, \c S. Raianu, \emph{Hopf Algebras: an introduction}. Vol. 235. Lecture Notes in Pure and Applied Math. Vol.235, Marcel Dekker, New York, 2001.

\bibitem[GTN]{NT1}
J. Gomez-Torrecillas, C. N\u ast\u asescu, \emph{Quasi-co-Frobenius coalgebras}, J. Algebra 174 (1995), 909-923.

\bibitem[GMN]{NT2}
J. Gomez-Torrecillas, C. Manu, C. N\u ast\u asescu, \emph{Quasi-co-Frobenius coalgebras II}, Comm. Algebra Vol 31, No. 10, pp. 5169-5177, 2003.

\bibitem[GNT]{NT3}
J. Gomez-Torrecillas, C. N\u ast\u asescu, B. Torrecillas, \emph{Localization in coalgebras. Applications to finiteness conditions}, eprint arXiv:math/0403248, http://arxiv.org/abs/math/0403248.

\bibitem[I]{I}
M.C.Iovanov, \emph{Co-Frobenius Coalgebras}, to appear, J. Algebra; eprint arXiv:math/0604251\\
http://xxx.lanl.gov/abs/math.QA/0604251.

\bibitem[I0]{I0}
M.C.Iovanov, \emph{Characterization of PF rings by the finite topology on duals of $R$ modules}, An. Univ. Bucure\c sti Mat. 52 (2003), no. 2, 189-200.

\bibitem[K1]{K1}
I. Kaplansky, \emph{Modules over Dedekind rings and valuation rings}, Trans. Amer. Math. Soc. 72 (1952) 327-340.

\bibitem[K2]{K2}
I. Kaplansky, \emph{A characterization of Pr$\rm\ddot{u}$fer domains}, J. Indian Math. Soc. 24 (1960) 279-281.

\bibitem[L]{L}
B.I-Peng Lin, \emph{Semiperfect coalgebras}, J. Algebra {\bf 30} (1974), 559-601.

\bibitem[McL]{McL} 
S. Mac Lane, \emph{Categories for the Working Mathematician}, Second Edition, Springer-Verlag, New York, 1971.

\bibitem[Mc1]{Mc1}
S. Mac Lane, \emph{Duality for groups}, Bull. Am. Math. Soc. {\bf 56}, 485-516 (1950).

\bibitem[NT]{NT}
C. N\u ast\u asescu, B. Torrecillas, \emph{The splitting problem for coalgebras}, J. Algebra {\bf 281} (2004), 144-149.

%\bibitem[NT0]{NT0}
%C. N\u ast\u asescu, B. Torrecillas, \emph{Symmetric Coalgebras}, J. Algebra {\bf 279} (2004), 326-344.

\bibitem[NTZ]{NTZ}
N\u ast\u asescu, B. Torrecillas, Y. Zhang, \emph{Hereditary Coalgebras}, Comm. Algebra {\bf 24} (1996), 1521-1528.

\bibitem[Rot]{Rot}
J. Rotman, \emph{A characterization of fields among integral domains}, An. Acad. Brasil Cienc. 32 (1960) 193-194.

\bibitem[T1]{T1}
M.L. Teply, \emph{The torsion submodule of a cyclic module splits off}, Canad. J. Math. XXIV (1972) 450-464.

\bibitem[T2]{T2}
M.L. Teply, \emph{A history of the progress on the singular splitting problem}, Universidad de Murcia, Departamento de ${\rm\acute{A}}$lgebra y Fundamentos, Murcia, 1984, 46pp.

\bibitem[T3]{T3}
M.L. Teply, \emph{Generalizations of the simple torsion class and the splitting properties}, Canad. J. Math. 27 (1975) 1056-1074.

\end{thebibliography}
\end{document}